\documentclass[onefignum,onetabnum,a4paper]{article}

\usepackage{amssymb}
\usepackage{amsmath}
\usepackage{epsfig}
\usepackage{psfrag}
\usepackage{subfigure}

\usepackage{xcolor}

\usepackage[a4paper, total={6in, 8in}]{geometry}

\newcommand{\fref}[1]{Fig.~\ref{#1}}

\newcommand{\Fref}[1]{Figure~\ref{#1}}

\title{Signatures consistent with multi-frequency tipping in the Atlantic meridional overturning circulation} 

\author{{\sc Andrew Keane\footnote{corresponding author:  {\tt andrew.keane@ucc.ie}}} \\ School of Mathematical Sciences and\\ Environmental Research Institute, University College Cork, \\College Road, Cork, Ireland
\and
{\sc Bernd Krauskopf} \\ Department of Mathematics, The University of Auckland, \\ Private Bag 92019, Auckland 1142, New Zealand
\and
{\sc Timothy M. Lenton} \\ Global Systems Institute, University of Exeter, \\Exeter, EX4 4QE, UK
}

\date{October 2020}

\begin{document}

\maketitle

%%%%%%%%%%%%%%%%%%%%%%%%%%%%%%%%%%%%%%%%%%%%%%%%%%%%%%%
\begin{abstract}
The early detection of tipping points, which describe a rapid departure from a stable state, is an important theoretical and practical challenge. Tipping points are most commonly associated with the disappearance of steady-state or periodic solutions at fold bifurcations. We discuss here multi-frequency tipping (M-tipping), which is tipping due to the disappearance of an attracting torus. M-tipping is a generic phenomenon in systems with at least two intrinsic or external frequencies that can interact and, hence, is relevant to a wide variety of systems of interest. We show that the more complicated sequence of bifurcations involved in M-tipping provides a possible consistent explanation for as yet unexplained behavior observed near tipping in climate models for the Atlantic meridional overturning circulation. More generally, this work provides a path towards identifying possible early-warning signs of tipping in multiple-frequency systems.
\end{abstract}
%%%%%%%%%%%%%%%%%%%%%%%%%%%%%%%%%%%%%%%%%%%%%%%%%%%%%%%

%%%%%%%%%%%%%%%%%%%%%%%%%%%%%%%%%%%%%%%%%%%%%%%%%%%%%%%
\paragraph*{Introduction}
\label{sec:intro}

A tipping point describes a sudden and often irreversible transition of a system from one state to another. Due to their potentially drastic impacts, tipping points have garnered much attention over recent years, particularly in climatology \cite{thompson13,dekker18,lenton19}, but also in biology and social science; see, for example, \cite{scheffer13,centola18}.

From a dynamical systems perspective, tipping may be induced by a bifurcation, the effect of noise or the rate of parameter change \cite{ashwin12}.
The ``classic'' tipping scenario involves a fold or saddle-node bifurcation of \emph{equilibria} \cite{thompson13} or, more recently, of \emph{periodic orbits} \cite{medeiros17,bathiany18}. Mathematically, the mechanism of tipping due to a fold is the same: a stable and a saddle steady state or periodic orbit merge and disappear at a single bifurcation point when a parameter changes slowly.

It is certainly the case that many systems of interest, in particular climate systems, are subject to forcing at various time scales and/or different types of feedback mechanisms; for example, see Ref.~\cite{stocker01} and references therein. This implies that the attractor undergoing a tipping event may well be more complex than an equilibrium or a periodic orbit (representing a single frequency). In particular, it is a natural and generic phenomenon that the dynamics evolves on an invariant \emph{torus} (representing two independent frequencies). Hence, in multi-frequency systems it is natural to study tipping involving tori. A much used approach is to average the dynamics to consider a simplified equilibrium case, as, for example, in zero-dimensional energy balance models for the overall temperature on the Earth's surface \cite{fraedrich79}. However, information is lost in this reduction and we will demonstrate that the actual dynamics on tori may offer important insights into tipping events in multi-frequency systems.  This is due to the fact that, generally, an attracting torus cannot simply lose stability and cease to exist at a single fold bifurcation point \cite{vitolo11}. Instead, the torus loses smoothness and ``breaks up'' in a complicated sequence of bifurcations \cite{baesens07,bakri14} before disappearing. We refer to this entire transition as multi-frequency tipping (M-tipping). 

Arguably, the most well-studied climate tipping scenario is that of a collapse of the Atlantic meridional overturning circulation (AMOC). The existence of bistable (``on'' and ``off'') solutions and a corresponding hysteresis loop dates back to the conceptual model of Stommel \cite{stommel61}, and it has since been observed in models across the complexity hierarchy. In the past, the transitions between these states have been described by fold bifurcations of equilibria. However, it is widely recognized that the AMOC exhibits internal modes of oscillatory variability at various frequencies --- most notably the Atlantic Multidecadal Oscillation (AMO) \cite{dijkstra06}. 

We show here that M-tipping is consistent with observations of possible precursors of tipping in intermediate complexity Earth-system models that have been used to investigate AMOC collapse, including a `step-down' weakening of AMOC before its total collapse. Hence, M-tipping may open the door to a much more detailed understanding of the dynamics leading to AMOC collapse and similar tipping scenarios in other fields of application.

%%%%%%%%%%%%%%%%%%%%%%%%%%%%%%%%%%%%%%%%%%%%%%%%%%%%%%%
\paragraph*{M-tipping in a conceptual DDE model}
\label{sec:dde}

In Ref.~\cite{keane18}, the authors studied the sudden disappearance of dynamics on a torus in a basic yet quite general model for the interaction of negative feedback with periodic forcing, which are both common ingredients in complex systems for generating oscillatory modes. In the context of AMOC it could be thought of as representing the interaction of any of the various negative feedback mechanisms \cite{buckley16,weijer19}, such as the temperature-advection feedback, with oscillatory behavior that could represent, for example, AMO or solar variability.

The basic conceptual model takes the form 
\begin{equation}
\label{eq:dde}
dh(t) = ( -\tanh[\kappa h(t-\tau_n)]+c \cos(2\pi t) )dt+\epsilon dW.
\end{equation}
The first term represents delayed negative feedback with delay time $\tau_n$ and nonlinearity strength $\kappa$, and the second term periodic forcing of strength $c$; together they form the delay differential equation (DDE) that was first introduced in Ref.~\cite{ghil08} and studied in Ref.~\cite{keane18}. The third term of additive white noise of strength $\epsilon$ has been included to demonstrate that signatures of M-tipping are robustly observed also in the presence of noise. Throughout, the noise level $\epsilon$ is taken to be small enough to ensure that the observed behavior is still driven chiefly by deterministic dynamics. In order to observe torus break-up, we consider $\tau_n$ and $c$ as bifurcation parameters and set $\kappa=11$, as in Refs.~\cite{ghil08,keane18}. Simulations in Ref.~\cite{ghil08} indicate that folding tori only occur for sufficiently large $\kappa$. Since this parameter determines a relative timescale of switching between two different levels of feedback, smaller values of $\kappa$ may also generate the phenomena discussed here if both $h$ and $t$ are rescaled accordingly.

Note that we do not attempt to calibrate the model to any one context. Rather, we use it here to demonstrate the significance of torus dynamics for AMOC tipping; in this context the variable $h(t)$ represents the strength of the Atlantic meridional overturning circulation.

%%%%%%%%%%%%%%%%%%%%%%%%%%%%%%%%%%%%%%%%%%%%%%%%%%%%%%%
\begin{figure}[t]
  \centering
\includegraphics[width=0.6\columnwidth]{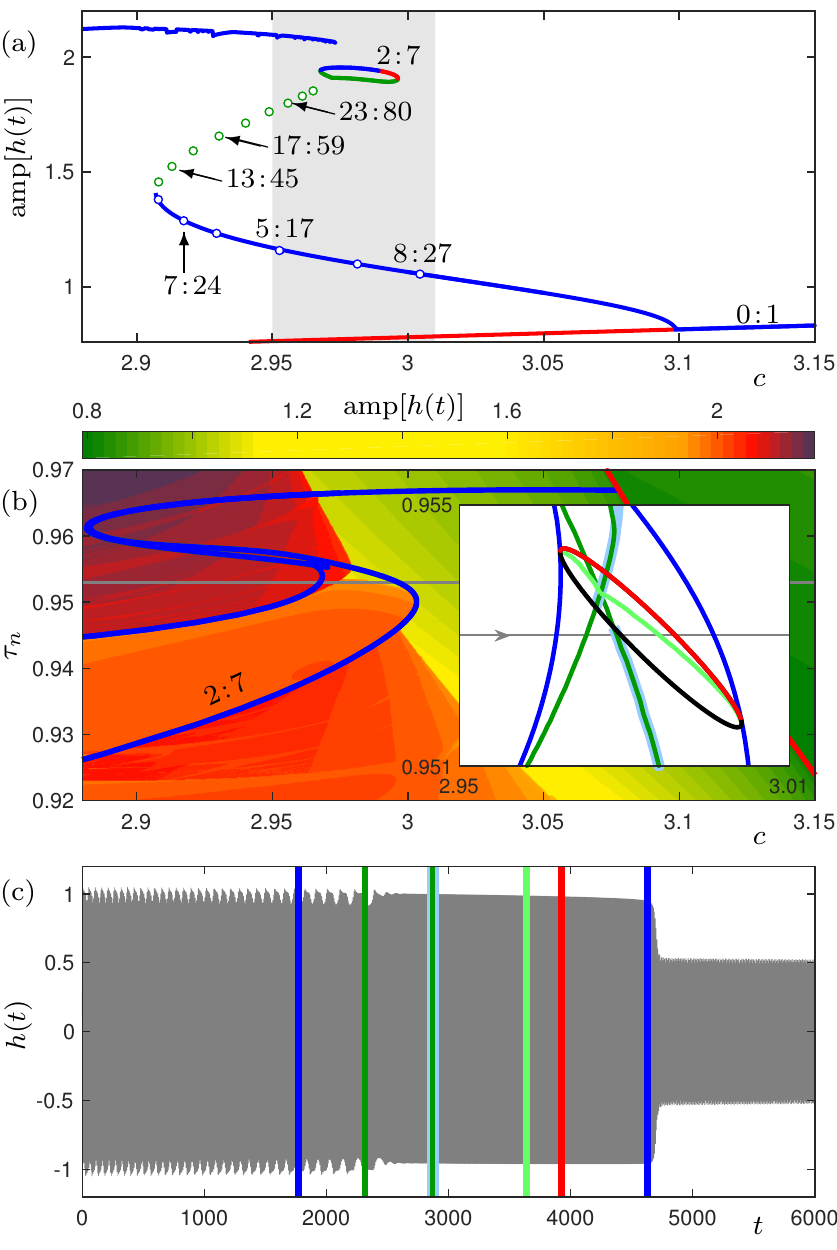}
  \caption{Dynamics near an M-tipping event of Eq.~\eqref{eq:dde} for $\epsilon=0$. (a) Hysteresis loop in parameter $c$ for $\tau_n=0.953$, represented by their amplitude $\text{amp}[h(t)]$. Stable solutions are shown in blue, while those with one and two unstable Floquet multipliers are shown in green and red, respectively. 
  (b) The ${2\!:\!7}$ resonance tongue in the $(c,\tau_n)$-plane. The color scheme shows amplitudes of stable solutions taken from simulations while slowly increasing $c$. The inset is an enlargement of the resonance tongue where it folds; see text for a list of bifurcation types.
  (c) Time series with linearly increasing $c$ corresponding to the gray-shaded region in panel~(a).
}
  \label{fig:gzt}
\end{figure}
%%%%%%%%%%%%%%%%%%%%%%%%%%%%%%%%%%%%%%%%%%%%%%%%%%%%%%%

We first briefly review torus break-up and M-tipping in the conceptual DDE model \eqref{eq:dde} in the absence of noise, that is, for $\epsilon=0$; see Ref.~\cite{keane18} for more details and, for example, Ref.~\cite{kuznetsov04} for background information on bifurcation theory. \Fref{fig:gzt}(a) depicts a hysteresis loop for $\tau_n=0.953$ of solutions on tori, which are represented by the amplitude of their time evolution $h(t)$. Along the blue curves one finds stable tori. The green and red curves represent periodic orbits with one and two unstable Floquet multipliers, respectively, found with the continuation software DDE-Biftool \cite{engelborghs00,sieber14}. The solutions on the tori are periodic when the two frequencies involved are locked into a rational ratio; some of which are labeled in \fref{fig:gzt}(a). Conversely, if the frequencies have an irrational ratio, the trajectories on the torus never close, and the torus is said to be \emph{quasi-periodic}. Most of the labeled periodic solutions exist across very small intervals of $c$ values and are represented by blue or green circles. The ${2\!:\!7}$ locked solutions, on the other hand, exist across a relatively large range of $c$ and appear in panel~(a) as an isola bounded by folds of periodic orbits. 
The green circles represent points along a branch of saddle tori that appears to connect the upper and lower branches of stable tori. We stress that this hysteresis loop is different from those for equilibria or periodic orbits. Namely, the connection between the different branches is not due to a simple fold but involves more complicated bifurcations and associated dynamics.

\Fref{fig:gzt}(b) shows the associated Arnold or resonance tongue in the $(c,\tau_n)$-plane to which the ${2\!:\!7}$ solutions belong. The background colors represent the amplitudes of simulated solutions for very slowly increasing values of $c$ across a range of fixed $\tau_n$ values. The gray line indicates the value of $\tau_n$ in panel~(a) and the inset details the dynamics within the resonance tongue across the shaded range of $c$-values in panel~(a). The curves shown here are: fold bifurcation of periodic orbits (dark blue), torus bifurcation (red), neutral saddle periodic orbit (black), homoclinic transition (light green), heteroclinic transition (dark green) and folding tori (light blue). To be precise, the homoclinic and heteroclinic curves represent transitions through tangles involving the stable and unstable manifolds of the ${2\!:\!7}$ locked saddle periodic orbits on the torus; since the first and last tangencies lie extremely close together, we represent each of these transitions by a single curve. The ${2\!:\!7}$ resonance tongue, bounded by two dark blue curves, folds twice in the $c$-direction. As mentioned earlier, a pair of stable and saddle tori cannot simply collide in the same way that equilibria or periodic orbits can in a fold bifurcation. The bifurcation structure in the inset (also called a Chenciner bubble) describes how the stable torus breaks up during such a folding and transitions to a saddle torus; see Ref.~\cite{keane18} for details. 
Notice in \fref{fig:gzt}(b) that, as $c$ increases, there is an associated sudden jump in the amplitude of the simulated solutions from large amplitude (red/orange) to small amplitude (yellow/green) behavior. 

The effect of torus break-up on the variable $h(t)$ is demonstrated in \fref{fig:gzt}(c) with a simulation that slowly traverses the gray-shaded range of $c$ in panel~(a), beginning on the upper branch of tori. Initially, one observes in panel~(c) two frequencies as the modulation of the amplitude. After about $t=2400$, the trajectory passes through a heteroclinic transition and approaches an attracting periodic orbit, which then loses stability near $t=4000$ at a subcritical torus bifurcation. Since the periodic orbit is only very weakly unstable, the trajectory diverges only very slowly from the periodic orbit until it completely disappears at the boundary of the ${2\!:\!7}$ resonance tongue near $t=4600$. The trajectory then approaches a remote attractor, namely, the lower branch of tori shown in \fref{fig:gzt}(a).

Generally, how observable the effects of torus break-up are depends on the size of the associated Chenciner bubble, which in turn depends on the frequency ratio ${p\!:\!q}$. For example, the effects of torus break-up are not so obvious on the left-hand side of the hysteresis loop in \fref{fig:gzt}(a) because these ${p\!:\!q}$ lockings are of a higher order $q$. Other factors include the parameter drift rate and effects of noise.

%%%%%%%%%%%%%%%%%%%%%%%%%%%%%%%%%%%%%%%%%%%%%%%%%%%%%%%
\paragraph{Signatures of M-tipping in AMOC}
\label{sec:amoc}

%%%%%%%%%%%%%%%%%%%%%%%%%%%%%%%%%%%%%%%%%%%%%%%%%%%%%%%
\begin{figure}[tbp]
  \centering
\includegraphics[width=0.6\columnwidth]{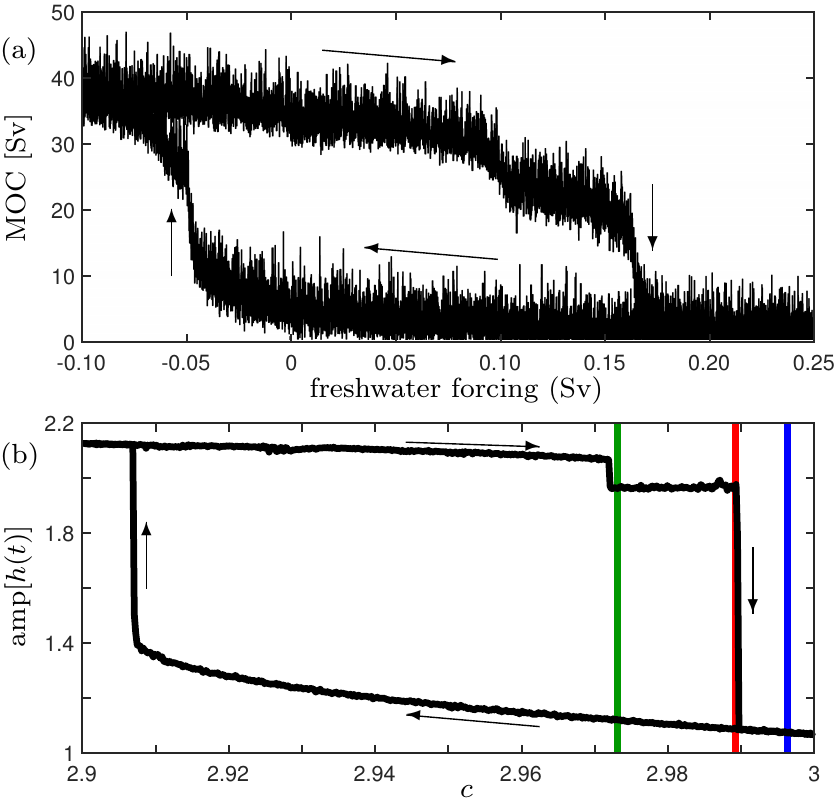}\\
  \caption{Hysteresis loops of (a) the GENIE-2 model as reproduced from Ref.~\cite{lenton09}, and (b) the conceptual DDE model \eqref{eq:dde} for $\tau_n=0.953$, $\Delta c=1\times10^{-6} \text{ per time unit}$ and $\epsilon=0.0005$.
}
  \label{fig:hyst}
\end{figure}
%%%%%%%%%%%%%%%%%%%%%%%%%%%%%%%%%%%%%%%%%%%%%%%%%%%%%%%

Hysteresis loops of AMOC have been demonstrated across the hierarchy of climate models. One such example from Ref.~\cite{lenton09}, based on the intermediate complexity GENIE-2 ocean-atmosphere model, is shown in \fref{fig:hyst}(a). It shows the maximum Atlantic meridional overturning circulation (MOC) as freshwater forcing is slowly varied according to the directions of the arrows. An interesting feature of the upper branch is described as a ``step slowdown'' in Ref.~\cite{lenton09}, whereby there is a relatively small, yet significant, drop in MOC as the forcing reaches approximately 0.10 Sv, before the full collapse at approximately 0.17 Sv (1 Sverdrup (Sv) $= 10^6 \text{m}^3\text{s}^{-1}$). In fact, a modest drop in MOC preceding a tipping event is a feature observed in several models \cite{ganopolski01,rahmstorf95,rahmstorf05}.

Since the maximum is always taken, the MOC observable may be interpreted as a measure of the size of an underlying oscillating attractor. This is supported by the fact that, if the maximum is not taken, but rather overturning circulation is considered for a \emph{fixed} depth and latitude, then variability on interdecadal time scales is observed in the output of a globally-coupled general circulation model (GCM) \cite{dong05}. 
To mimic the MOC observable in a simple way, we show in \fref{fig:hyst}(b) the observable $\text{amp}[h(t)]$ of Eq.~\eqref{eq:dde} (with nonzero noise level $\epsilon$), which now approximates the amplitude of the underlying attractor as max minus min across a window of 200 time units.
As $c$ increases, before the big drop there is a small drop in $\text{amp}[h(t)]$ where the trajectory passes a heteroclinic transition (green line). This happens shortly before the bifurcation point is reached due to anticipation by noise.  Similarly, the big drop at the end of the hysteresis loop already occurs near the subcritical torus bifurcation (red line). Larger values of $\epsilon$ result in the transitions occurring for slightly smaller values of $c$.
Clearly, the dynamics associated with M-tipping, demonstrated in \fref{fig:hyst}(b), offer a possible mechanism for the ``step slowdown'' observed in panel~(a). In M-tipping, the observation that the small drop in MOC only occurs on the upper branch of the hysteresis loop implies that the underlying Chenciner bubble, where the torus breaks up, is much smaller on the lower branch, so that its effects are not observed. This is indeed the case for Eq.~\eqref{eq:dde} as \fref{fig:hyst}(b) and \fref{fig:gzt}(a) show.

%%%%%%%%%%%%%%%%%%%%%%%%%%%%%%%%%%%%%%%%%%%%%%%%%%%%%%%
\begin{figure}[tbp]
  \centering
\includegraphics[width=0.6\columnwidth]{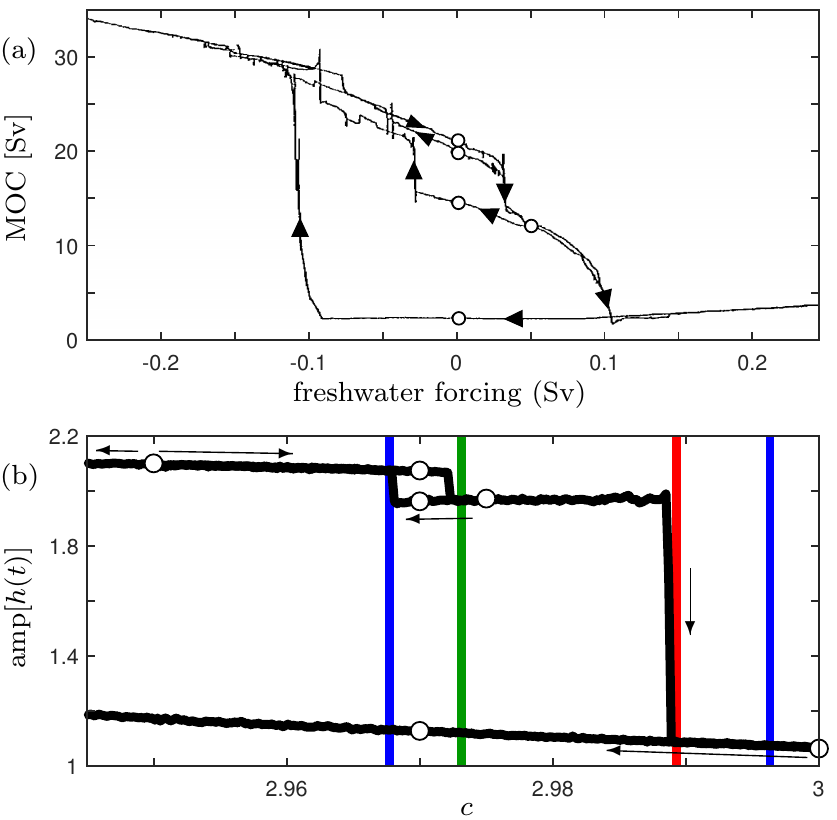}
  \caption{Additional multistability near the tipping event observed in (a) a global ocean GCM as reproduced from Ref.~\cite{rahmstorf95}, and (b) the conceptual DDE model \eqref{eq:dde} for $\tau_n=0.953$, $\Delta c=1\times10^{-6} \text{ per time unit}$ and $\epsilon=0.0005$. 
}
  \label{fig:tri}
\end{figure}
%%%%%%%%%%%%%%%%%%%%%%%%%%%%%%%%%%%%%%%%%%%%%%%%%%%%%%%

\Fref{fig:tri}(a) shows a similar freshwater forcing experiment by Rahmstorf \cite{rahmstorf95} with a coarse-resolution, yet highly-sophisticated, ocean general circulation model (GCM); here the arrows indicate the direction in which the forcing parameter is changed and the open circles represent stable solutions after very long transient times. Apart from a similar small drop in MOC (referred to as an ``interesting discontinuity''), Rahmstorf observed a region along the upper branch with two distinctly different ``on'' solutions with MOC values of about 20 and 15 Sv and one ``off'' solution near 0 Sv.
This behavior is reproduced for the conceptual DDE model \eqref{eq:dde} and the observable $\text{amp}[h(t)]$ in \fref{fig:tri}(b).
%%, where open circles are again stable solutions after very long transients. 
We also observe a range of tristability, here bounded by a fold bifurcation of periodic orbits (blue line) and a heteroclinic transition (green line), where three simultaneously stable solutions are indicated by open circles. In this $c$-range, there is locally a coexistence of a stable torus and a stable periodic orbit solution as two ``on'' solutions, showing that
M-tipping also provides a possible mechanism for the additional multistability observed near the main tipping event in the ocean GCM.

In Ref.~\cite{rahmstorf95} Rahmstorf identified that the small drop in MOC results from a change in convection patterns, representing the local shutdown of deep water formation in the Labrador Sea. This process was simulated in detail during freshwater forcing experiments with the intermediate complexity ECBilt-CLIO model in Ref.~\cite{schulz07}. 

%%%%%%%%%%%%%%%%%%%%%%%%%%%%%%%%%%%%%%%%%%%%%%%%%%%%%%%
\begin{figure}[tbp]
  \centering
\includegraphics[width=0.6\columnwidth]{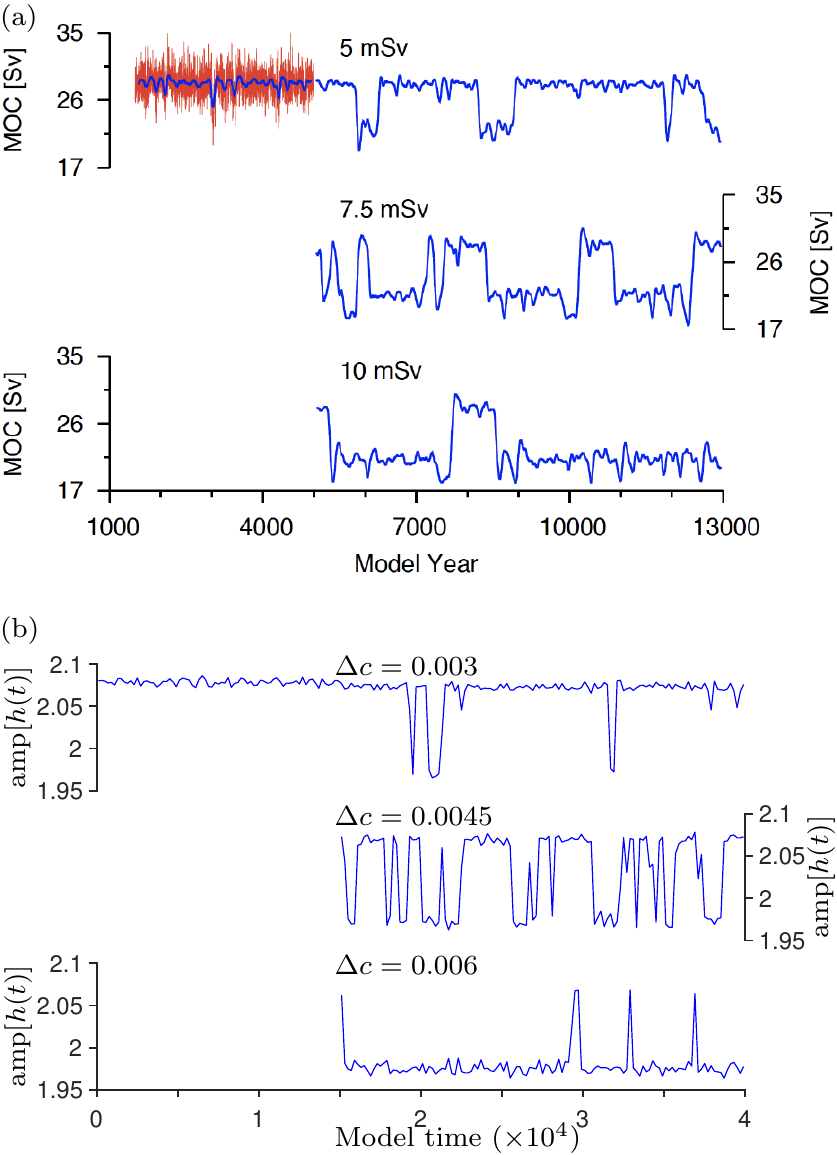}
  \caption{Intermittent times series in (a) ECBilt-CLIO model with a freshwater perturbation of 5 mSv, 7.5 mSv and 10 mSv at 5000 years as reproduced from Ref.~\cite{schulz07}, and (b) the conceptual DDE model \eqref{eq:dde} for $c=2.966,\tau_n=0.953$ and $\epsilon=0.001$, with a perturbation $\Delta c$ in parameter $c$ of 0.003, 0.0045 and 0.006.
}
  \label{fig:intermittent}
\end{figure}
%%%%%%%%%%%%%%%%%%%%%%%%%%%%%%%%%%%%%%%%%%%%%%%%%%%%%%%

\Fref{fig:intermittent}(a) from Ref.~\cite{schulz07} shows responses in MOC to three different freshwater perturbations. The first 5000 years are simulated for the unperturbed case, showing unsmoothed data (red) and data smoothed with a 101-year Hanning filter (blue). For each perturbation, of 5 mSv, 7.5 mSv and 10 mSv, the response is bi-modal and switches between a higher and lower state that correspond to Labrador Sea convection being ``on'' and ``off', respectively. This more gradual approach towards the small drop in MOC reveals an intermittent transition: the trajectory gradually spends more and more time in the lower state, until the transition is complete.

\Fref{fig:intermittent}(b) shows the corresponding transition in Eq.~\eqref{eq:dde} near the heteroclinic transition.
After $1.5\times10^4$ unperturbed time units at $c=2.966$, for increasing $\Delta c$ perturbations of 0.003, 0.0045 and 0.006
one clearly observes an intermittent transition due to noise-induced jumps between a stable torus coexisting with a stable periodic orbit.
As the dynamics on the torus becomes increasingly slow in regions of phase space where it is close to the periodic orbit, jumps from the torus to the periodic orbit are increasingly favored. 
Different values of $\epsilon$ effect the range of $c$ values for which the intermittent transition occurs, but not the occurence of the phenomenon itself.
We conclude that M-tipping is also consistent with the intermittent transition observed in Ref.~\cite{schulz07} when convection in the Labrador Sea shuts down.

%%%%%%%%%%%%%%%%%%%%%%%%%%%%%%%%%%%%%%%%%%%%%%%%%%%%%%%
\paragraph*{Discussion}
\label{sec:discuss}

The dynamics associated with M-tipping provide an elegant explanation for phenomena observed in three different models of AMOC. Our working hypothesis is that the two convection zones in the Nordic Seas and the Labrador Sea each contribute towards the formation of feedback loops, similar to those represented by simple ocean box-models, such as in Ref.~\cite{stommel61}. Two feedback loops can naturally give rise to two-frequency dynamics on a torus, which then must break up and pass through a heteroclinic transition. In the process, there is a switch to one-frequency dynamics, which appears to be achieved in AMOC by the shutdown of convection in one of its convection zones, namely, in the Labrador Sea. Future work on AMOC will focus on modeling specific feedback loops and the analysis of tori in models of intermediate complexity.

As basic as it may seem, the conceptual DDE model considered here highlights the potentially crucial role of torus dynamics for understanding certain tipping phenomena. As such we hope that it will provide motivation for further studies that embrace the multi-frequency nature of complex systems. Additional bifurcation structure associated with M-tipping, such as the heteroclinic transition central to the results presented here, could provide new early-warning indicators of an approaching tipping event. Whether such indicators can be identified reliably in different scenarios requires further study.

\paragraph*{Acknowledgements}
The authors thank Henk Dijkstra and Jan Sieber for helpful discussions. The research of BK was supported by Royal Society Te Ap\={a}rangi Marsden Fund grant \#19-UOA-223. TL thanks the University of Exeter for teaching leave, the University of Auckland for hosting him, and the Leverhulme Trust for funding (RPG-2018-046).

\bibliographystyle{siam}
\bibliography{kkl_mft_refs}

\end{document}